\documentclass[12pt]{amsart} 

\usepackage[latin9]{inputenc}
\usepackage{amssymb}
\usepackage{latexsym}
\usepackage{amsmath}
\usepackage{amsthm}
\usepackage{amsfonts}

\usepackage{tikz}
\usepackage{dsfont}
\usepackage{mathrsfs}
\usepackage{bbm}
\usepackage{color}

\allowdisplaybreaks

\newtheorem{theorem}{Theorem}[section]

\newtheorem*{lemma*}{Lemma}
\newtheorem{corollary}[theorem]{Corollary}
\newtheorem{proposition}[theorem]{Proposition}
\newtheorem*{proposition*}{Proposition}

\theoremstyle{definition}

\newtheorem*{acknowledgements*}{Acknowledgements}
\theoremstyle{remark}
\newtheorem{remark}[theorem]{Remark}

\numberwithin{equation}{section}

\newcommand{\N}{\mathbb N}
\newcommand{\C}{\mathbb C}

\newcommand{\D}{\mathbb D}

\newcommand{\DD}{\overline{\mathbb D}}

\newcommand{\cno}{\C^{\N_0}}
\newcommand{\no}{\N_0}
\newcommand{\ce}{\mathcal{C}}

\newcommand{\ces}{Ces\`{a}ro }
\newcommand{\subin}{_{n=0}^\infty}
\newcommand{\subi}{_{n=0}^\infty}

\newcommand{\dpp}{d_{pp}}
\newcommand{\dpe}{d_p}
\newcommand{\mdp}{\mathscr{M}(d_p)}
\newcommand{\bodp}{\mathscr{A}(S,\dpe)}
\newcommand{\ldp}{\mathscr{L}(d_p)}
\newcommand{\hdp}{H(d_p)}
\newcommand{\mhdp}{\mathscr{M}(H(d_p))}
\newcommand{\hmdp}{H(\mathscr{M}(d_p))}
\newcommand{\modp}{\mathscr{M}_{\mathrm{op}}(d_p)}
\newcommand{\hinf}{H^\infty(\D)}
\newcommand{\lwp}{\ell^1(w_p)}
\newcommand{\hd}{H(\D)}
\newcommand{\hdd}{H(\DD)}
\newcommand{\mdone}{\mathscr{M}(d_1)}
\newcommand{\mlp}{\mathscr{M}(\ell^p)}

\title
[Convolution  in  dual Ces\`{a}ro  sequence spaces]
{Convolution in dual Ces\`{a}ro   sequence spaces}

\author[G.~P.  Curbera]{Guillermo P. Curbera}
\address{Facultad de Matem\'aticas \& IMUS,
Universidad de Sevilla, 
Calle Tarfia s/n,  Sevilla 41012, Spain}
\email{curbera@us.es}

\author[W.~J. Ricker]{Werner J. Ricker}
\address{Math.--Geogr. Fakult\"at, Katholische Universit\"at
Eichst\"att--Ingolstadt, D--85072 Eichst\"att, Germany}
\email{werner.ricker@ku.de}

\thanks{Both authors acknowledge the support of the 
Mathematisches Forschungsinstitut Oberwolfach via the
Research in Pairs Program (March, 2022). 
The first author also acknowledges the support  of 
PID2021-124332NB-C21 FEDER/Ministerio de Ciencia e Innovaci\'on and 
FQM-262 (Spain).
}

\date{\today}
\subjclass[2010]{Primary 47B37, 47L10; Secondary  46B45, 47A10.}
\keywords{Banach algebra, convolution, dual Ces\`{a}ro sequence space, multiplier, spectrum.}

\begin{document}

\begin{abstract}
We investigate  convolution operators in the sequence spaces
$d_p$, for $1\le p<\infty$. These  spaces, for $p>1$, arise as dual spaces of the \ces sequence
spaces $ces_p$ thoroughly investigated by G.~Bennett. A detailed study is also made
of the algebra of those sequences which convolve $d_p$ into $d_p$. It turns out
that such multiplier spaces exhibit features which are very different to
the classical multiplier spaces of $\ell^p$.
\end{abstract} 

\maketitle


\section{Introduction}
\label{S1}


In 1966, in a celebrated paper, \cite{nikolskii}, N.~K.~Nikolskii  
initiated the study of multipliers acting on the classical sequence spaces 
$\ell^p=\ell^p(\no)$, with $\no=\{0,1,2,\dots\}$, where 
\begin{equation*}
\ell^p:=\Big\{a=(a_n)\subi\in\cno:\sum\subin|a_k|^p<\infty\Big\},
\quad 1\le p<\infty. 
\end{equation*}
A sequence $b=(b_n)\subi\in\cno$ defines a \textit{multiplier} on $\ell^p$ 
if the \textit{convolution}  $a*b\in\cno$, defined by
\begin{equation}\label{eq-2-10}
(a*b)_n:=\sum_{j=0}^na_jb_{n-j},\quad n\in\no,
\end{equation}
belongs to $\ell^p$, for every $a\in\ell^p$. 
The \textit{multiplier algebra} $\mlp$ of $\ell^p$ is the collection of all such $b\in\cno$.
Nikolskii established the following fundamental properties of  these multiplier algebras:
\begin{itemize}
\item[a)] $\ell^1\subsetneq \mlp\subsetneq\ell^p$, for $1<p<\infty$;
\item[b)] $\mlp=\mathscr{M}(\ell^{p'})$, for $1/p+1/p'=1$;  
\item[c)] $\mathscr{M}(\ell^{p_1})\subsetneq \mathscr{M}(\ell^{p_2})$, for $1\le p_1<p_2\le2$.
\end{itemize}
These multiplier algebras, except when $p\in\{1,2\}$, are not well 
understood and their  investigation is  far from finalized. 
Important contributions were made by Vinogradov, Verbitskii 
and others; see, for example,   \cite[\S6.41--6.43]{bottcher-silbermann}, and
\cite{cheng-etal} for a recent account of the state of the art.


The \ces sequence spaces $ces_p$, for $1<p<\infty$, are intimately connected to the spaces
$\ell^p$ via the \ces averaging operator which maps each element of $\ell^p$ to the sequence of its averages (again an element of $\ell^p$). The spaces $ces_p$
were throughly investigated by G.~Bennett, \cite{bennett}; see also
\cite{grosse-erdmann} and the references therein. They have the property that 
$\ell^p\subsetneq ces_p$, for all $1<p<\infty$. However, in contrast to  $\ell^p$,
the situation regarding the multipliers of $ces_p$ is completely
different: the multiplier algebra $\mathscr{M}(ces_p)=\ell^1$,
for \textit{every} $1<p<\infty$, \cite[Theorem 4.1]{curbera-ricker-ieot-2014}.

The purpose of this note is to investigate the multiplier algebras
$\mdp$ of the sequence spaces $\dpe$, also spaces 
closely related to $\ell^p$, which are defined by
\begin{equation}\label{eq-1-1}
d_p:=\Big\{a=(a_n)\subi\in\cno:\sum\subin\sup_{k\ge n}|a_k|^p<\infty\Big\},
\quad 1\le p<\infty.
\end{equation}
They were defined and studied by G.~Bennett, \cite{bennett},  
when he obtained a tractable identification of the dual Banach space of $ces_p$. 
More precisely,  the dual Banach space $(ces_p)^*$ is isomorphic to $d_q$, 
for $p\in(1,\infty)$, where $\frac1p+\frac1q=1$; \cite[Corollary 12.17]{bennett}.
Despite having  similarities in their definition, the spaces
$\ell^p$ and $\dpe$ are rather different. A significant difference is that the canonical vectors
$e_n:=(\delta_{n,k})_{k=0}^\infty$, for $n\in\no$, are all
\textit{unit vectors} in every space $\ell^p$, for $p\in[1,\infty]$, but they have norm
$\|e_n\|_{\dpe}=(n+1)^{1/p}$ whenever $1\le p<\infty$ and $n\in\no$.
For further properties of the spaces $\dpe$, see \cite{bonet-ricker},
for example. Note that $\dpe\subsetneq \ell^p \subsetneq ces_p$,
for $1<p<\infty$.

The multiplier algebras $\mdp$ of  $\dpe$ consist of all $b\in\cno$
which convolve $\dpe$ into itself. Differences between the spaces $\ell^p$ and $\dpe$
induce  drastically different features between their respective multiplier spaces $\mlp$
and $\mdp$. In contrast to property a) above, we have that
$$
\mdp\subsetneq\ell^1\quad\mathrm{and}\quad 
\mathscr{M}(d_1)=d_1\subsetneq \mdp
\subsetneq\dpe,\quad1<p<\infty;
$$
see Theorem \ref{t-3-5} and Corollary \ref{c-3-6}. 
That is, \textit{all} the spaces $\mdp$ are inside $\ell^1$.
In contrast to properties b) and c) above, it turns out that 
$$
\mathscr{M}(d_{p_1})\subsetneq \mathscr{M}(d_{p_2}),
\quad 1\le p_1<p_2<\infty;
$$
see Theorem \ref{t-3-8}. That is, there is no largest space with the role
that $\mathscr{M}(\ell^{2})$ has in the $\ell^p$ setting.

As for $\mlp$, with  $p\not\in\{1,2\}$, no characterization of the entire
algebra $\mdp$ is known (except for $p=1$). Nevertheless, we
devote some effort to
identify natural classes of elements which do belong to $\mdp$.
For example, the weighted Banach algebra $\lwp$ with
$w_p(n)=(n+1)^{1/p}$ for $n\in\no$ is contained
in $\mdp$ for every $1\le p<\infty$; see Proposition \ref{p-3-7}.
A characterization of those elements from $\ell^1$ which belong to $\mdp$ is presented in
Theorem \ref{t-4-1}. 
A more tractable sufficient condition for a sequence $b\in\ell^1$
to be a multiplier for $\dpe$, in terms of its coefficients, namely that
$$
\sum_{n=0}^\infty 2^{np}  \sup_{2^{n}\le k< 2^{n+1}}|b_k|^p<\infty,
$$
is established in Theorem \ref{t-4-2}.

Together with $\mdp$ we also consider the associated algebra
$\modp$ of all (necessarily) bounded, linear convolution operators $T_b$ 
on $\dpe$ induced by the elements $b$ of $\mdp$; see Section \ref{S2}
for the definitions. As for the spaces $\ell^p$, the right-shift operator $S$ (which maps
an element $(a_0,a_1,\dots)$ to $(0,a_0,a_1,\dots)$ ) also plays
an important role for the spaces $\dpe$. For instance, it turns out that the commutant algebra
$\modp^c$ of $\mdp$ equals
\begin{equation}\label{eqn-1-5}
\modp^c=\Big\{T\in\ldp: TS=ST\Big\},\quad 1\le p<\infty,
\end{equation}
where $\ldp$ is the space of all bounded linear operators of $\dpe$
into itself. 
A crucial difference between the $\ell^p$ and the $\dpe$ setting is that
the operator norm of $S^n\in\ldp$ equals $(n+1)^{1/p}$
for each $n\in\no$ and $1\le p<\infty$, whereas $S^n\in\mlp$
is an isometry  for all such $n$ and $p$. 
Consequences of \eqref{eqn-1-5} are that 
$\mdp$ is complete for the weak operator topology (cf.\ Section \ref{S3}) and that 
the spectrum of an operator
in the unital, commutative Banach algebra $\modp$, for $1\le p<\infty$,
coincides with its spectrum as an element of $\ldp$. The
topic of the spectrum of operators belonging to $\modp$ is pursued
in the final section. Of particular relevance are the distinct subspaces
$d_1, \lwp$ and $\dpp\cap\ell^1$ of $\mdp$ because,
if $b=(b_n)\subi$ belongs to any one of these subspaces, then
the corresponding multiplier operator $T_b\in\modp$ can be
approximated in the operator norm by the polynomial operators
$\{\sum_{k=0}^nb_kS^k\}\subi$; see 
Remark \ref{r-6-6}(ii)  and Proposition \ref{p-5-7}.

The paper is organized as follows. Section \ref{S2} presents the necessary preliminaries
required in the sequel.
Section \ref{S3} treats various relevant properties of the operator algebras $\modp$, whereas
Section \ref{S4} concentrates on the multiplier algebras $\mdp$. In 
Section \ref{S5} we identify various subspaces of  $\mdp$.
The final Section \ref{S6} is devoted to 
spectral and Banach algebra properties of $\modp$.


\section{Preliminaries}
\label{S2}


For each $p\in[1,\infty)$ the sequence space $\dpe$ defined in \eqref{eq-1-1}
is a  Banach space for   the norm
\begin{equation}\label{eq-2-1}
\|a\|_{d_p}:=\Big(\sum\subin \sup_{k\ge n}|a_k|^p\Big)^{1/p},\quad a\in d_p.
\end{equation}
A direct consequence of \eqref{eq-2-1} is that $d_p\subseteq \ell^p$ with a continuous inclusion.
Given $a=(a_n)\subi\in\ell^\infty$, the \textit{least decreasing majorant of $a$} is the sequence
$\hat a:=(\sup_{k\ge n}|a_k|)\subi$,  \cite[(3.7)]{bennett}. 
Then, $a\in\dpe$ precisely when $\hat a\in\ell^p$ and $\|a\|_{\dpe}=\|\hat a\|_{p}$,
where $\|\cdot \|_{p}$ is the usual norm in $\ell^p$. 
The canonical vectors  $\{e_n:n\in\no\}$ satisfy
\begin{equation*}
\|e_n\|_{\dpe}=\|\widehat{e_n}\|_{\ell^p}=\big\|(1,\dots,1,
\overbrace{1}^{\text{position } n},
0,0,\dots)\big\|_{p}=(n+1)^{1/p}.
\end{equation*}
For every $p\in[1,\infty)$, the vectors $\{e_n:n\in\no\}$ form an unconditional basis in $\dpe$, 
\cite[Proposition 2.1]{bonet-ricker};  see Section \ref{S4} for the case $p=1$.

A combination of Cauchy's condensation test for series and Abel's summation formula 
implies the following two  useful equivalent expressions for the norm \eqref{eq-2-1} in $d_p$: 
\begin{align}\label{eq-2-3}
\|a\|_{d_p} &
\asymp \left(\sup_{k\ge0} |a_k|^p+\sum\subin 2^n \sup_{2^n\le k<2^{n+1}} |a_k|^p\right)^{1/p},
\\ \|a\|_{d_p}  & \asymp \label{eq-2-4}
\left(\sup_{k\ge0} |a_k|^p+\sup_{k\ge1} |a_k|^p+ \sum_{n=0}^\infty 2^n \sup_{2^n< k\le2^{n+1}} |a_k|^p\right)^{1/p},
\end{align}
where $A\asymp B$ means that there exist absolute constants $c,C>0$ 
such that $cA\le B\le CA$;
see also \cite[Example 13.2]{grosse-erdmann} and \cite[(3)]{belinskii}.

As noted in Section \ref{S1}, the space  $d_q$ is isomorphic to  $(ces_p)^*$, where $ces_p$,  \cite{bennett},  is defined, for each $1< p\le\infty$, by
\begin{equation}\label{eq-2-5}
ces_p:=\Big\{a=(a_n)\subi\in\cno:
\|a\|_{ces_p}:=\Big(\sum\subi\Big(\frac{1}{n+1}\sum_{k=0}^{n}|a_k|\Big)^p\Big)^{1/p}\Big\},
\end{equation}
that is, $a\in ces_p$ if and only if $\big(\frac{1}{n+1}\sum_{k=0}^{n}|a_k|\big)\subi\in\ell^p$.


The \textit{convolution} of $a,b\in\cno$ is the sequence $a*b\in\cno$ defined by
\eqref{eq-2-10}. According to Section \ref{S1} 
the \textit{multiplier algebra} 
\begin{equation*}\label{eq-2-11}
\mdp:=\Big\{b\in\cno:a*b\in\dpe, \forall a\in\dpe\Big\}.
\end{equation*}
Each $b\in\mdp$ defines a
convolution operator  $a\mapsto a*b\in\dpe$, for $a\in\dpe$, which is continuous
(due to the closed graph theorem). The multiplier
algebra $\mdp$  endowed with the  norm 
\begin{equation}\label{eq-2-7}
\|b\|_{\mdp}:= \sup_{0\not=a\in\dpe}\frac{\|a*b\|_{\dpe}}{\|a\|_{\dpe}},
\end{equation}
is a Banach algebra; see Section \ref{S3}. Since $e_0\in\dpe$ satisfies $e_0*b=b$
for every $b\in\cno$, it is clear that
$\mdp\subseteq\dpe$. 
This implies (as mentioned above) that $\mdp$  is a unital, commutative algebra under convolution.
Moreover, for each $b\in\mdp$, we have that
$\|b\|_{\dpe}=\|e_0*b\|_{\dpe}/ \|e_0\|_{\dpe}\le\|b\|_{\mdp}$.
Since $\|e_0\|_{\mdp}=1=\|e_0\|_{\dpe}$, it follows that the operator norm of the
natural inclusion $\mdp\subseteq\dpe$ is precisely 1.


\section{The  operator algebra $\modp$}
\label{S3}


Convolution operators on $\dpe$ will be considered within the 
unital (non-commutative) Banach algebra $\ldp$  
of all bounded linear operators on $\dpe$ equipped with the operator norm. 
Given $b\in\mdp$, denote by $T_b$   the convolution operator defined 
by $T_b(a):=a*b\in\dpe$, for each $a\in\dpe$, and set
\begin{equation*}
\modp:=\Big\{T_b\in\ldp:b\in\mdp\Big\}.
\end{equation*}
Observe that $\|T_b\|_{\modp}=\|b\|_{\mdp}$ for all $b\in\mdp$. 
Clearly, $\modp$ is a commutative, 
unital subalgebra of $\ldp$,
with the identity operator $I=T_{e_0}$ as its unit. Equipped with
the operator norm from $\ldp$, which we denote by 
$\|\cdot\|_{\modp}$, it becomes a normed algebra.

The \textit{commutant algebra}  of $\mdp$ is defined by 
\begin{equation*}
\modp^c:=\Big\{R\in\ldp: T_bR=RT_b, \;\forall b\in\mdp\Big\}. 
\end{equation*}


The right-shift  $S\colon\dpe\to\dpe$ is the linear map given by
\begin{equation*}
Sa=(0,a_0,a_1,\dots)=e_1*a=T_{e_1}a,\quad a\in\dpe.
\end{equation*}
It follows, for $n\in\no$, that 
\begin{equation*}
S^na=(0,\dots,0,\overbrace{a_0}^{\text{position } n},a_1,\dots)=e_n*a=T_{e_n}a,\quad a\in\dpe.
\end{equation*}
Direct calculation yields $\|e_n\|_{\dpe}=\|S^n\|_{\modp}=(n+1)^{1/p}$,
for $n\in\no$ and $p\in[1,\infty)$; see \cite[Lemma 4.12]{curbera-ricker-jmaa-2022}.
This is distinctly different to the situation for the spaces $\ell^p$, where
$\|e_n\|_{p}=\|S^n\|_{\mathscr{L}(\ell^p)}=1$,
for all $n\in\no$ and $p\in[1,\infty]$.


\begin{proposition}\label{p-3-1}
Let $p\in[1,\infty)$. Then
\begin{equation}\label{eq-3-1}
\modp=\Big\{R\in\ldp: RS=SR\Big\}.
\end{equation}
Moreover,
\begin{equation}\label{eq-3-2}
\modp=\modp^c=\modp^{cc}.
\end{equation}
\end{proposition}

\begin{proof}
Let $T\in\ldp$ satisfy $TS=ST$ and set $b:=Te_0\in\dpe$. 
Since $e_1=Se_0$, we have $Te_1=TSe_0=STe_0=Sb=b*e_1$. In a similar way, using
$e_{n+1}=Se_n$, it follows that $Te_n=S^nb=b*e_n$ for all $n\in\no$. Hence,
$Ta=b*a$ for all $a$ belonging to the linear span of
$\{e_n:n\in\no\}$. Since the  canonical vectors $\{e_n:n\in\no\}$ form a basis for 
$\dpe$, for every $a=(a_n)\subi\in\dpe$ we have $a^N\to a$ in 
$\dpe$, where $a^N=\sum_{j=0}^N a_je_j$. Then $Ta^N\to Ta$ in $\dpe$ and so 
$b*a^N\to Ta$ in $\dpe$. Since convergence in $\dpe$ implies coordinatewise
convergence, for each fixed $n\in\no$, we have
$$
(b*a^N)_n=\Big(b*\sum_{j=0}^N a_je_j\Big)_n=
\Big(\sum_{j=0}^N a_j(b*e_j)\Big)_n\to (Ta)_n.
$$
Note, for $N\ge n$, that 
$$
\Big(\sum_{j=0}^N a_j(b*e_j)\Big)_n= \Big(\sum_{j=0}^n a_j(b*e_j)\Big)_n
= \Big(\sum_{j=0}^n a_jS^jb\Big)_n=(b*a)_n.
$$
Hence, $(b*a)_n=(Ta)_n$ for $n\in\no$, that is, $b*a=Ta$ and so, $b*a\in\dpe$.
Since $a\in\dpe$ is arbitrary, we have $b\in\mdp$ and $T=T_b$.

The reverse inclusion in \eqref{eq-3-1} follows easily as $S=T_{e_1}\in\modp$.

Since $\modp$ is commutative, it is contained in $\modp^c$. On the other hand,
if $R\in\modp^c$, then $S=T_{e_1}$ implies that $RS=SR$ and so, by \eqref{eq-3-1},
the operator $R\in\modp$. Hence, $\modp=\modp^c$. It then follows that
$$
\modp^{cc}=(\modp^c)^c=\modp^c=\modp,
$$
which is precisely \eqref{eq-3-2}.
\end{proof}


\begin{remark}\label{r-3-2}
(i) For the spaces $\ell^p$ in place of $\dpe$, with $p\in[1,\infty)$, the
identity \eqref{eq-3-1} is known, \cite[Theorem 2(2)]{nikolskii}. Also,
for $ces_p$ in place of $\dpe$, with $p\in(1,\infty)$, the same proof
as in Proposition \ref{p-3-1} applies to show that identities \eqref{eq-3-1} and \eqref{eq-3-2}
hold. However, unlike for $\ell^p$ and $\dpe$, we have the remarkable fact that
$$
\mathscr{M}_{\text{op}}(ces_p)=\big\{T_b:b\in\ell^1\big\},\quad p\in(1,\infty),
$$
and that $\|T_b\|_{ces_p\to ces_p}=\|b\|_1$ for $a\in\ell^1$;
see \cite[Theorem 4.1]{curbera-ricker-ieot-2014}.

(ii) In view of \eqref{eq-3-2} it is well known that $\modp$ is
\textit{inverse closed} in $\ldp$, \cite[I Proposition 2.3]{bonsall-duncan},
that is, if $T\in\modp$ is invertible in $\ldp$, then its inverse operator $T^{-1}\in\ldp$
actually belongs to $\modp$. In particular, the spectrum  $\sigma(R;\modp)$ 
of an operator $R\in\modp$ coincides with its spectrum $\sigma(R;\ldp)$ as an element of $\ldp$. For the definition of the spectrum of an element in a unital
Banach algebra we refer to \cite{bonsall-duncan}, \cite{naimark}, for example.
\end{remark}


\begin{corollary}\label{c-3-3}
For each $p\in[1,\infty)$ the algebra $\modp$ is closed in $\ldp$ for the weak
operator topology and hence, also for the strong operator topology and the 
operator norm
topology. In particular, $\modp$ is a commutative  Banach algebra (i.e., it is complete).
\end{corollary}

\begin{proof}
Let $\{T^{(\alpha)}\}\subseteq\modp$ be a net and $T\in\ldp$ such that 
$T^{(\alpha)}\overset{\alpha}{\to} T$ for the weak operator topology. Proposition \ref{p-3-1}
yields $T^{(\alpha)}S=ST^{(\alpha)}$ for all $\alpha$. Fix $a\in\dpe$ and $y^*\in \dpe^*$.
Then, with $S^*\in\mathscr{L}(\dpe^*)$ denoting the adjoint operator of $S$, we have
\begin{align*}
\langle STa,y^*\rangle&=\langle Ta,S^*y^*\rangle=\lim_\alpha\langle T^{(\alpha)}a,S^*y^*\rangle
\\ & =\lim_\alpha\langle ST^{(\alpha)}a,y^*\rangle
=\lim_\alpha\langle T^{(\alpha)}Sa,y^*\rangle=\langle TSa,y^*\rangle.
\end{align*}
It follows that $TS=ST$ and hence, that $T\in\modp$.
\end{proof}


\section{The  multiplier algebra $\mdp$}
\label{S4}


In this section we study various properties of 
the multiplier algebras $\mdp$. We begin with 
$p=1$   which is simpler and is already known.  Recall that
\begin{equation*}
d_1:=\Big\{a=(a_n)\subi\in\cno:\|a\|_{d_1}:=\sum_{n=0}^\infty\sup_{k\ge n}|a_k|<\infty\Big\},
\end{equation*}
which can be traced back to the 
work of Beurling, \cite{beurling}; see Remark \ref{r-4-1} below.
The canonical vectors  $\{e_n:n\in\no\}$ form an unconditional 
basis in $d_1$. This follows 
from a necessary condition for  a sequence $a=(a_n)\subi$ to belong to $d_p$,
namely that
$$
\lim_{n} n\sup_{k\ge n} |a_k|^p=0,
$$
which is a consequence of Pringsheim's theorem for convergent series  
of positive decreasing terms. Indeed, for $a\in d_1$, we have for $N\to\infty$ that
\begin{align*}
\Big\|a-\sum_{n=0}^N a_ne_n\Big\|_{d_1} &=
\Big\|\Big(\sup_{k\ge N+1} |a_k|, \dots, 
\overbrace{\sup_{k\ge N+1} |a_k|}^{\text{position } N+1}
, \sup_{k\ge N+2} |a_k|,\dots\Big\|_{\ell^1}
\\ & =N\sup_{k\ge N+1} |a_k|+\sum_{n=N+1}^\infty\sup_{k\ge n} |a_k| \to 0.
\end{align*}
The bounded multiplier test ensures the unconditionality of the basis. 
The space $d_1$ is known to be
an algebra for convolution with unit $e_0$ 
(see the proof of \cite[Proposition 1]{belinskii}). So, 
$\mdone$ and $d_1$ coincide as sets
and have equivalent norms, that is, for some $C>0$ we have
\begin{equation*}\label{eqn-2-2}
\|b\|_{d_1}\le \|b\|_{\mdone} \le C\,\|b\|_{d_1},\quad b\in d_1,
\end{equation*}
where we have used $\|b\|_{d_1}=\|T_be_0\|_{d_1}$ and \eqref{eq-2-7}.
In particular, $\mdone\subsetneq\ell^1$ (since 
$|a|\le\hat{a}$ and \cite[Remark 4.20(i)]{curbera-ricker-jmaa-2022} 
imply that $d_1\subsetneq\ell^1$).


\begin{remark}\label{r-4-1}
A result of Beurling concerning the absolute convergence of contracted Fourier
series is based on imposing  on the Fourier coefficients $(a_n)_{-\infty}^\infty$ 
of an integrable function on $[0,2\pi]$ the   condition
$$
\sum\subin \sup_{|k|\ge n} |a_k|<\infty,
$$
\cite[Theorem V]{beurling}. Note that $d_1$ corresponds to this condition  
when $a_n=0$ for $n<0$.
\end{remark}

The following result already indicates how different the
multiplier algebras $\mdp$ and $\mlp$  are.


\begin{theorem}\label{t-3-5}
For each   $p\in[1,\infty)$, the following continuous inclusion holds:
$$
\mdp\subseteq\ell^1.
$$
\end{theorem}

\begin{proof}
For $p=1$ this is  $d_1\simeq\mdone\subseteq\ell^1$. For $p\in(1,\infty)$,
let $0\not=b\in\mdp$. Denote by  $n_0$ the smallest $n\in\no$ such that $b_n\not=0$. 
Fix $n\ge n_0$. For any $a\in\dpe$, it follows from \eqref{eq-2-3} that
\begin{align*}
\|a*b\|_{\dpe}^p & \ge 2^n \sup_{2^n\le k<2^{n+1}} |(a*b)_k|^p
\ge 2^n |(a*b)_{2^n}|^p = 2^n  \Big|\sum_{j=0}^{2^n}b_j a_{2^n-j}\Big|^p.
\end{align*}
Define $a=(a_n)\subi\in\dpe$ via $a_{2^n-j}=|b_j|/b_j$ for $0\le j\le 2^n$ (with 
$a_{2^n-j}=0$ if $b_j=0$) and $a_j=0$ for $j>2^n$. Then
$$
\sum_{j=0}^{2^n}b_j a_{2^n-j}=\sum_{j=0}^{2^n} |b_j|.
$$
Note that $\|a\|^p_{\dpe}\le(2^n+1)$.
Consequently,
\begin{align*}
\|b\|_{\mdp}^p = \sup_{0\not=a\in\dpe} \frac{\|a*b\|_{\dpe}^p}{\|a\|_{\dpe}^p} \ge 
\frac{\displaystyle 2^n  \Big(\sum_{j=0}^{2^n}|b_j|\Big)^p}{2^n+1} \ge
\frac12  \Big(\sum_{j=0}^{2^n}|b_j|\Big)^p.
\end{align*}
It follows that $b\in\ell^1$ and $ \sum_{j=0}^{\infty}|b_j| \le 2^{1/p} \|b\|_{\mdp}$. 
\end{proof}


\begin{corollary}\label{c-3-6}
Let $p\in(1,\infty)$. The following assertions hold.
\begin{itemize}
\item[(i)] $\mdp\subsetneq \dpe$.
\item[(ii)] $\mdp\not=\ell^1$. 
\end{itemize}
\end{corollary}

\begin{proof}
(i) We have already seen in Section \ref{S2} that  $\mdp\subseteq\dpe$. 
Let $a=(1/(n+1))\subi$. Since it is a decreasing sequence and $a\in\ell^p$, we see that $a\in\dpe$.
However, since $a\not\in\ell^1$, we have $a\not\in\mdp$.
Note that $a$ is the sequence of Taylor coefficients of the analytic function
$\log(1-z)\not\in\hinf$.

(ii) Suppose that $\mdp=\ell^1$.
Since  $\mdp\subseteq \dpe$ this would imply that $\ell^1\subseteq \dpe$, 
which is not the case; see \cite[Remark 2.8(i)]{bonet-ricker}.
\end{proof}


Consider the weight $w_p:=((n+1)^{1/p})\subi$ and the 
corresponding weighted $\ell^1$-space
\begin{equation*}
\ell^1(w_p):=\Big\{(a_n)\subi: \sum\subi (n+1)^{1/p}|a_n|<\infty\Big\},
\end{equation*}
equipped with the norm $\|a\|_{1,w_p}:=\sum\subi (n+1)^{1/p}|a_n|$.
Observe that $w_p(m+n)\le w_p(m)w_p(n)$ for all $m,n\in\no$.


\begin{proposition}\label{p-3-7}
For each $p\in[1,\infty)$ the following continuous embedding holds:
\begin{equation*}
\ell^1(w_p)\subseteq \mdp.
\end{equation*}
\end{proposition}

\begin{proof}
Let  $m\in\no$. The  canonical vector $e_m\in\dpe$ defines a multiplier in $d_p$. 
Indeed, fix $a\in\dpe$. Since
$$
e_m\ast a=(\overbrace{0,\dots,0}^m,a_0,a_1,\dots), 
$$
the least decreasing majorant of $e_m\ast a$ is
$$
(e_m\ast a)\,\hat{}=\bigg(
\overbrace{\sup_{k\ge0}|a_k|, \dots,\sup_{k\ge0}|a_k|}^{m+1}, \sup_{k\ge1}|a_k|, \dots\bigg).
$$
But,  $a\in d_p$ and so $\hat{a}\in \ell^p$. By the previous identity it is 
clear that $(e_m\ast a)\,\hat{}\in \ell^p$ and
\begin{equation*}
\big\|e_m\ast a\big\|_{d_p}=
\Big\|(e_m\ast a)\,\hat{}\ \Big\|_p= \Big(m\big(\sup_{k\ge0}|a_k|\big)^p
+\|a\|_{d_p}^p\Big)^{1/p}.
\end{equation*}
In particular,  $\big\|e_m\ast a\big\|_{d_p} \le (m+1)^{1/p} \|a\|_{d_p}$.
Consequently,  $e_m\in\mdp$ and $\|e_m\|_{\mdp}\le (m+1)^{1/p}$. 
This bound  is sharp as
can be seen by selecting $a=e_0$, in which case 
$e_m*e_0=e_m$ with $\hat{e_m}=\sum_{n=0}^{m}e_n$.
So, $\|e_m\|_{\mdp}\ge(m+1)^{1/p}$. Hence, $\|e_m\|_{\mdp}=(m+1)^{1/p}$.

Let $a=(a_n)\subi\in\ell^1(w_p)$. Consider in $\mdp$ the series $\sum\subi a_ne_n$. 
It is absolutely convergent in $\mdp$ because
$$
\sum\subi\|a_ne_n\|_{\mdp}=\sum\subi|a_n|\|e_n\|_{\mdp}
=\sum\subi|a_n|(n+1)^{1/p}=\|a\|_{1,w_p}.
$$
Since the space $\mdp\simeq\modp$ is complete (cf.\  Corollary \ref{c-3-3}), it follows that 
the series is convergent in $\mdp$.
\end{proof}


\begin{theorem}\label{t-3-8}
Let $1\le p_1<p_2<\infty$. Then $\mathscr{M}(d_{p_1})
\subsetneq\mathscr{M}(d_{p_2})$.
In particular, $d_1\subseteq\mdp$ for all $1\le p<\infty$.
\end{theorem}

\begin{proof}
We first show, for $1\le p_1<p_2<\infty$, that  $d_{p_2}$ is an interpolation
space between $d_{p_1}$ and $\ell^\infty$. More precisely, we will show that 
\begin{equation}\label{eq-new}
(d_{p_1})^\theta (\ell^\infty)^{1-\theta}=d_{p_2},
\quad\mathrm{for}\quad \theta:=\frac{p_1}{p_2}\in(0,1),
\end{equation}
where $(d_{p_1})^\theta (\ell^\infty)^{1-\theta}$ is 
a Calder\'on space, \cite[13.5]{calderon}. 
Observe that each space $\dpe$  is the Tandori 
space corresponding to $\ell^p$ since, in the notation of
\cite{lesnik-maligranda}, for $a=(a_n)\subi\in\ell^\infty$, 
we have $\widetilde a =\hat a$,  \cite[\S1]{lesnik-maligranda}.  Recall that
$\hat a$ is the decreasing majorant of $a$ (cf.\ \S2). 
Consequently,  $\widetilde{\ell^p}=\dpe$, for $1\le p<\infty$;
see  \cite[(1.6)]{lesnik-maligranda}. It is clear  that
$\widetilde{\ell^\infty}=\ell^\infty$. 

Theorem 4 in \cite{lesnik-maligranda} states, for suitable spaces $X_0, X_1$ and
an adequate  function $\varphi$
(cf.\ \cite[\S3]{lesnik-maligranda}), that 
$$
\varphi(\widetilde{X_0}, \widetilde{X_1})=[\varphi(X_0, X_1)]\;\widetilde{}.
$$
We apply this result to the spaces $X_0=\ell^{p_1}$, $X_1=\ell^{\infty}$ and 
the function $\varphi(s,t):=s^\theta t^{1-\theta}$ with 
$\theta:=p_1/p_2\in(0,1)$. Then,
$\widetilde{X_0}=d_{p_1}$, $\widetilde{X_1}=\ell^{\infty}$
and $\varphi(X_0, X_1)=(\ell^{p_1})^\theta(\ell^\infty)^{1-\theta}=\ell^{p_2}$, so that
$[\varphi(X_0, X_1)]\;\widetilde{}=d_{p_2}$. Thus, the equality \eqref{eq-new} follows.

Let $b\in\mathscr{M}(d_{p_1})$. Then $T_b\colon d_{p_1}\to d_{p_1}$. 
Theorem \ref{t-3-5} yields that 
$b\in \ell^1$. This implies, for $a\in\ell^\infty$ and every $n\in\no$, that
$|(a*b)_n|\le\sum_{j=0}^n |a_jb_{n-j}|\le \|a\|_\infty\|b\|_1$, that is, $T_ba\in\ell^\infty$. 
Hence, $T_b\colon \ell^\infty\to \ell^\infty$. 
The equality \eqref{eq-new} implies that $d_{p_2}$ is a Calder\'on $\theta$-space
for $d_{p_1}$ and $\ell^\infty$. So, $d_{p_2}$ is an interpolation space
between $d_{p_1}$ and $\ell^\infty$, 
\cite[33.5]{calderon}. This yields that 
$T_b\colon d_{p_2}\to d_{p_2}$, that is, $b\in\mathscr{M}(d_{p_2})$.


To show that $\mathscr{M}(d_{p_1})\not=\mathscr{M}(d_{p_2})$,
let $b=(b_n)\subi$ be defined by $b_n=2^{-k/p_1}$ when $n=2^k$ 
(for $k\in\no$) and $b_n=0$ 
otherwise. Since $\frac{1}{p_1}>\frac{1}{p_2}$, it follows that
$$
\sum\subi |b_n|(n+1)^{1/p_2}=\sum_{k=0}^\infty \frac{(2^k+1)^{1/p_2}}{2^{k/p_1}}<\infty,
$$
and so $b\in\ell^1(w_{p_2})$. From Proposition \ref{p-3-7} we have 
$\ell^1(w_{p_2})\subseteq\mathscr{M}(d_{p_2})$, that is,
$b\in\mathscr{M}(d_{p_2})$. 
However, $b\not\in d_{p_1}$ because
$$
\sum\subi 2^n \sup_{2^n\le k<2^{n+1}} |b_k|^{p_1}=
\sum\subi 2^n  |b_{2^n}|^{p_1}=\sum\subi \frac{2^n}{(2^{n/p_1})^{p_1}}=\infty.
$$
Since $\mathscr{M}(d_{p_1})\subseteq d_{p_1}$, it 
follows that $b\not\in\mathscr{M}(d_{p_1})$. Hence, 
$\mathscr{M}(d_{p_1})\subsetneq\mathscr{M}(d_{p_2})$.

By the discussion prior to Remark \ref{r-4-1} we have that $d_1=\mdone$, which implies
that $d_1\subseteq\mdp$ for all $1\le p<\infty$.
\end{proof}


\begin{remark}
(i) We  also refer to \cite[\S15 p.176]{maligranda} for spaces of the form $X_0^\theta X_1^{1-\theta}$
and \cite[Theorem 3]{shestakov} for an interpolation theorem for these spaces.

(ii) In the proof of Theorem \ref{t-3-8}, an alternative way of showing  that $d_{p_2}$ 
is an interpolation space between
$d_{p_1}$ and $\ell^\infty$, for $1\le p_1<p_2<\infty$, is via 
an interpolation result for Wiener-Beurling spaces.
More precisely,  Theorem 5.1(i) in \cite{nursultanov}
applied to $WB^{1/p_1}_{\infty,p_1}(\no)=d_{p_1}$, 
$WB^{0}_{\infty,\infty}(\no)=\ell^\infty$ and $WB^{1/p_2}_{\infty,p_2}
(\no)=d_{p_2}$ yields $(d_{p_1},\ell^\infty)_{1-\frac{p_1}{p_2},p_2}=d_{p_2}$.
\end{remark}


Let $\hd$ denote the space of all analytic functions on $\D$.
Consider the space of those functions in $\hd$ whose Taylor 
coefficients belong to   $\dpe$,
namely,  
\begin{equation*}\label{eq-2-8}
\hdp:=\Big\{f_a(z):=\sum\subi a_n z^n: (a_n)\subi\in\dpe\Big\}
\subseteq \hd,
\end{equation*}
where the notation $f_a$ indicates that $a=(a_n)\subi$ is the 
sequence of Taylor coefficients of $f_a$. Since $\dpe\subseteq\ell^\infty$, it is 
clear that $f_a$ is indeed analytic
in $\D$ for each $a\in\dpe$. 
The norm in $\hdp$ is defined by
\begin{equation*}\label{eq-2-9}
\|f_a\|_{\hdp}=\Big\|\sum\subi a_n z^n\Big\|_{\hdp}:=\|(a_n)\subi\|_{\dpe},\quad f_a\in\hdp.
\end{equation*}
Accordingly, as Banach spaces $\dpe$ and $\hdp$ are linearly isomorphic and isometric
via the map $a\leftrightarrow f_a$.
Consequently, the dual space $\hdp^*$ of $\hdp$ is isomorphic 
to the space $H(ces_q)$ of analytic functions
with Taylor coefficients in $ces_q$.

Given $z\in\D$ the point evaluation functional $\delta_{z}$ on $\hdp$,
for $p\in[1,\infty)$, is defined by
$$
f_a\in\hdp\longmapsto \delta_z(f_a):=f_a(z)=\sum\subi a_nz^n\in\C.
$$

\begin{proposition}\label{p-3-4} 
Let $p\in[1,\infty)$. For each $z\in\D$ the functional $\delta_z$ 
on $\hdp$ is linear and bounded,
that is, $\delta_z\in\hdp^*$. For $p\in(1,\infty)$ its norm satisfies
$$
\frac{1/p}{1-|z|}\Big(\sum\subi\bigg(\frac{1-|z|^{n+1}}{n+1}\bigg)^q\Big)^{1/q}
\le \|\delta_{z}\|_{\hdp^*}\le
\frac{(q-1)^{1/q}}{1-|z|}\Big(\sum\subi\bigg(\frac{1-|z|^{n+1}}{n+1}\bigg)^q\Big)^{1/q}, 
$$
where $\frac1p+\frac1q=1$. In particular,
$$
\frac1p\zeta(q)^{1/q} \le \|\delta_{z}\|_{\hdp^*}\le \frac{(q-1)^{1/q}}{1-|z|} \zeta(q)^{1/q}.
$$
For $p=1$, the functional  $\delta_z$ acting on $H(d_1)$ has norm one.
\end{proposition}


\begin{proof}
Fix $z\in\D$. Consider  $f_a(z)=\sum_{n=0}^\infty a_nz^n\in\hdp$. 
Then
\begin{equation}\label{eq-3-3}
\delta_{z}(f_a)=f_a(z)=\sum_{n=0}^\infty a_nz^n=
\big\langle\big(z^n\big)_{n=0}^\infty,(a_n)_{n=0}^\infty\big\rangle.
\end{equation}
For $p\in(1,\infty)$, we have $a\in\dpe$ and $(z^n)_{n=1}^\infty\in \ell^q\subseteq ces_q$, 
which is isomorphic to $\dpe^*$. 
Thus, $\delta_{z}$ acting on $\hdp$ can be identified with the sequence
$(z^n)_{n=0}^\infty\in (\dpe)^*$ acting
on $d_p$. Since $\hdp$ and $d_p$ are
isometric,  the norms of $\delta_{z}$  as an element of
$\hdp^*$ and of
$(z^n)_{n=0}^\infty$ as an element of $d_p^*$ coincide.
The equivalence of the norms between  $d_q$ and $(ces_p)^*$ is given by
\begin{equation}\label{eq-2-6}
\frac1q \|a\|_{d_q}\le \|a\|_{(ces_p)^*}\le (p-1)^{1/p} \|a\|_{d_q},\quad a\in (ces_p)^*,
\end{equation}
where $p$ and $q$ are conjugate indices, i.e.,  $\frac1p+\frac1q=1$,
\cite[p. 61 and Corollary 12.17]{bennett}. 
From \eqref{eq-2-6} it follows that  the 
equivalence  of the norms between $(\dpe)^*$ and $ces_q$ is given by
\begin{equation*}
\frac1p \|a\|_{ces_q}\le \|a\|_{(d_p)^*}\le (q-1)^{1/q} \|a\|_{ces_q},\quad a\in (d_p)^*.
\end{equation*}
In our case this yields
\begin{equation}\label{eq-3-4}
\frac1p\|(z^n)_{n=0}^\infty\|_{ces_q}
\le \|\delta_{z}\|_{\hdp^*}\le (q-1)^{1/q}\|(z^n)_{n=0}^\infty\|_{ces_q}.
\end{equation}
The norm of  $(z^n)_{n=0}^\infty$ in $ces_q$ is given by
\begin{align*}
\|(z^n)_{n=0}^\infty\|_{ces_q}^q&=
\sum\subi\bigg(\frac{1}{n+1}\sum_{k=0}^n|z^k|\bigg)^q
=\frac{1}{(1-|z|)^q}\sum\subi\bigg(\frac{1-|z|^{n+1}}{n+1}\bigg)^q .
\end{align*}
Since
$$
(1-|z|)^q\sum\subi\frac{1}{(n+1)^q}\le
\sum\subi\bigg(\frac{1-|z|^{n+1}}{n+1}\bigg)^q \le 
\sum\subi\frac{1}{(n+1)^q},
$$
we can conclude that
$$
\zeta(q)\le \|(z^n)_{n=0}^\infty\|_{ces_q}^q
\le \frac{\zeta(q)}{(1-|z|)^q}.
$$
The claim now follows from \eqref{eq-3-4}.

For $p=1$, from \eqref{eq-3-3} we have $a\in d_1$ and 
$(z^n)_{n=1}^\infty\in ces_\infty$, 
which is isometric to $d_1^*$, \cite[Remark 6.3]{curbera-ricker-ieot-2014}. 
Thus, $\delta_{z}$ acting on $H(d_1)$  can be identified with the sequence
$(z^n)_{n=0}^\infty\in (d_1)^*$ acting
on $d_1$. Hence, the norm of $\delta_{z}$ equals the norm of
$(z^n)_{n=0}^\infty$ in $ces_\infty$, that is,
$$
\big\|(z^n)_{n=0}^\infty\|_{ces_\infty}=\sup_{n\ge0}
\frac{1}{n+1}\sum_{k=0}^n|z|^k=1.
$$
\end{proof}

In view of the proof of the above result and the isomorphism $\dpe\simeq\hdp$,
it is clear, for each $z\in\D$, that $\delta_z\in\hdp^*$ corresponds to
the element of $\dpe^*$ given by $a\mapsto \sum\subi a_nz^n$, for $a\in\dpe$.


The Taylor coefficients of the pointwise product of two 
analytic functions $f_a$ and  $f_b$
in $\D$ are obtained via the convolution of $a$ and $b$, 
that is, $f_a f_b=f_{a*b}$. Consequently, the space 
\begin{equation*}
\mhdp:=\Big\{\varphi\in \hd: \varphi f\in\hdp, \forall f\in\hdp\Big\}
\end{equation*}
of analytic multipliers for $\hdp$ is linearly isomorphic and isometric to the space
$\hmdp$ of analytic functions on $\D$ with Taylor coefficients in the algebra 
$\mdp$, that is, to the algebra
\begin{equation*}
\hmdp:=\Big\{\varphi_a(z)=\sum_{n=0}^\infty a_n z^n: (a_n)\subi\in\mdp\Big\}
\subseteq \hd
\end{equation*}
equipped with the norm $\|\varphi_a\|_{\hmdp}:=\|a\|_{\mdp}$. 
Note the identification between 
$\mhdp$ and $\hmdp$.
Observe that $\hmdp\subseteq\hdp$ because $\mdp\subseteq\dpe$.

With obvious notation (that is, interchanging $d_p\leftrightarrow \ell^p$) 
it is known that
\begin{equation}\label{eq-3-5}
\ell^1\subseteq\mathscr{M}(\ell^p)\simeq\mathscr{M}(H(\ell^p))
\subseteq\hinf,\quad 1<p<\infty,
\end{equation}
where $\hinf$ is the space of all bounded analytic functions on $\D$, 
\cite[Theorem 4]{nikolskii}.
The containment in the right-side of \eqref{eq-3-5} 
can be sharpened when we consider 
$d_p$ in place of $\ell^p$. This is because $f_a(z)=\sum\subi a_nz^n\in\hmdp$ implies, via Theorem \ref{t-3-5}, that $a=(a_n)\subi\in\ell^1$, 
and so in \eqref{eq-3-5} we can replace
the space $\hinf$ by the  classical (one-sided) analytic 
Wiener algebra, \cite[\S11.6]{naimark}, denoted by $\ell^1_A$ in \cite{nikolskii}, consisting of all analytic functions 
on $\D$ with absolutely convergent Taylor coefficients. That is,
$$
d_1\subseteq\mdp\simeq \mhdp\subseteq\ell^1_A,\quad 1<p<\infty.
$$


\section{Subspaces of $\mdp$}
\label{S5}


Theorem \ref{t-3-5} shows for $b\in\cno$ that  a necessary condition for
being a multiplier  for $\dpe$ is that $b\in\ell^1$. This fact allows the formulation of a 
necessary and sufficient condition for $b\in\ell^1$ to belong to 
$\mdp$, which has the advantage 
that, for each $n\in\no$, in the $n$-th term
of the series in \eqref{eq-4-1}    below only the terms 
$b_j$ for $2^{n-1}<j< 2^{n+1}$ occur.


\begin{theorem}\label{t-4-1} 
Let $p\in(1,\infty)$ and $b\in\ell^1$. Then 
$b\in\mdp$ if and only if
\begin{equation}\label{eq-4-1}
\sum_{n=0}^\infty 2^n \sup_{2^n\le k<2^{n+1}} 
\Big|\sum_{\frac{k}{2}<j\le k}b_ja_{k-j}\Big|^p<\infty,
\quad a\in\dpe.
\end{equation}
\end{theorem}


\begin{proof}
Recall that $b\in\mdp$ if and only if 
$a*b\in\dpe$, for every $a\in\dpe$. This is equivalent, via \eqref{eq-2-3}, to
\begin{equation*}\label{eqn-2}
\sup_{n\ge0}|(a*b)_n|^p +
\sum_{n=0}^\infty 2^n \sup_{2^n\le k<2^{n+1}} \big|(a*b)_k\big|^p<\infty,
\quad a\in\dpe.
\end{equation*}
Since $b\in \ell^1$, given any $a\in d_p\subseteq\ell^p$ it follows 
that $a*b\in\ell^p$ and so, $a*b$ is bounded. 
Hence,   $b\in\mdp$ if and only if 
\begin{equation}\label{eq-4-2}
\sum_{n=0}^\infty 2^n \sup_{2^n\le k<2^{n+1}} \big|(a*b)_k\big|^p<\infty,
\quad a\in\dpe.
\end{equation}

First assume that the condition \eqref{eq-4-1} is satisfied.
To prove that  $b\in\mdp$ it suffices to establish \eqref{eq-4-2}. 
Let $a\in\dpe$. Then, for each $k\in\no$, we have
\begin{align}\label{eq-4-3}\nonumber
\big|(a*b)_k\big|=\Big|\sum_{j=0}^k b_ja_{k-j}\Big|
&=
\Big|\sum_{0\le j\le \frac{k}{2}} b_ja_{k-j}+\sum_{\frac{k}{2}<j\le k}b_ja_{k-j}\Big|  
\\ & \le 
\Big(\sum_{0\le j\le \frac{k}{2}} |b_j|\Big)\sup_{0\le j\le \frac{k}{2}} |a_{k-j}|
+\Big|\sum_{\frac{k}{2}<j\le k}b_ja_{k-j}\Big|  
\\ & \le \nonumber
\|b\|_1\sup_{\frac{k}{2}\le j\le k} |a_{j}|
+\Big|\sum_{\frac{k}{2}<j\le k}b_ja_{k-j}\Big| .
\end{align}
Fix $n\in\no$.  It follows from \eqref{eq-4-3} that
\begin{align}\label{eq-4-4}\nonumber
\sup_{2^n\le k<2^{n+1}} \big|(a*b)_k\big|^p 
&= \Big(\sup_{2^n\le k<2^{n+1}} \big|(a*b)_k\big|\Big)^p
\\ & \le 
\Big(\sup_{2^n\le k<2^{n+1}}\|b\|_1\sup_{\frac{k}{2}\le j\le k} |a_{j}|
+\sup_{2^n\le k<2^{n+1}}\Big|\sum_{\frac{k}{2}<j\le k}b_ja_{k-j}\Big|\Big)^p 
\\ & = \nonumber
\Big(\|b\|_1 \sup_{2^{n-1}\le k<2^{n+1}}|a_{j}|
+\sup_{2^n\le k<2^{n+1}}\Big|\sum_{\frac{k}{2}<j\le k}b_ja_{k-j}\Big|\Big)^p .
\end{align}
The  inequality \eqref{eq-4-4} implies that
\begin{align*}
\sum_{n=0}^\infty 2^n \sup_{2^n\le k<2^{n+1}} \big|(a*b)_k\big|^p
& \le
\sum_{n=0}^\infty 2^n   \Big(\|b\|_1 \sup_{2^{n-1}\le k<2^{n+1}}|a_{j}|
\\ &\quad +\sup_{2^n\le k<2^{n+1}}\Big|\sum_{\frac{k}{2}<j\le k}b_ja_{k-j}\Big|\Big)^p
\end{align*}
Applying  Minkowski's inequality yields 
\begin{align}\label{eq-4-5}\nonumber
\Big(\sum_{n=0}^\infty 2^n  
\sup_{2^n\le k<2^{n+1}}\big|(a*b)_k\big|^p\Big)^{1/p}
& \le
\Big(\sum_{n=0}^\infty 2^n   \big(\|b\|_1 \sup_{2^{n-1}\le k<2^{n+1}}|a_{j}|\big)^p\Big)^{1/p}
\\ &  
\quad +
\Big(\sum_{n=0}^\infty 2^n \sup_{2^n\le k<2^{n+1}}
\Big|\sum_{\frac{k}{2}<j\le k}b_ja_{k-j}\Big|^p\Big)^{1/p} .
\end{align}

The second term in the right-side of \eqref{eq-4-5}
is finite because of \eqref{eq-4-1}. 
Regarding the first term in the right-side of \eqref{eq-4-5}, note that
\begin{align}\label{eq-4-6}\nonumber
\sum_{n=0}^\infty 2^n   \sup_{2^{n-1}\le k<2^{n+1}}|a_{j}|^p
& \le\nonumber
\sum_{n=0}^\infty 2^n   \sup_{2^{n-1}\le k<2^{n}}|a_{j}|^p
+\sum_{n=0}^\infty 2^n   \sup_{2^{n}\le k<2^{n+1}}|a_{j}|^p
\\ & =
2\sum_{n=0}^\infty 2^{n-1}   \sup_{2^{n-1}\le k<2^{n}}|a_{j}|^p
+\sum_{n=0}^\infty 2^n   \sup_{2^{n}\le k<2^{n+1}}|a_{j}|^p
\\ & \le\nonumber
3\sum_{n=0}^\infty 2^n   \sup_{2^{n}\le k<2^{n+1}}|a_{j}|^p .
\end{align}
Then
\begin{equation}\label{eq-4-7}
\Big(\sum_{n=0}^\infty 2^n   \big(\|b\|_1 \sup_{2^{n-1}\le k<2^{n+1}}|a_{j}|\big)^p\Big)^{1/p}
\le \|b\|_1 3^{1/p}\Big(\sum_{n=0}^\infty 2^n   \sup_{2^{n}\le k<2^{n+1}}|a_{j}|^p\Big)^{1/p},
\end{equation}
which is also finite since $b\in\ell^1$ and $a\in\dpe$. 
Hence, \eqref{eq-4-2} is finite for every $a\in\dpe$ and so,   $b\in \mdp$.


Conversely,  we need to show that condition \eqref{eq-4-1} is necessary.
So, assume that  $b\in\mdp$.
Fix $a\in\dpe$. Then

\begin{align*}
\sum_{n=0}^\infty 2^n \sup_{2^n\le k<2^{n+1}} \Big|\sum_{\frac{k}{2}<j\le k}b_ja_{k-j}\Big|^p
& = 
\sum_{n=0}^\infty 2^n \sup_{2^n\le k<2^{n+1}} 
\Big|\sum_{0\le j\le k}b_ja_{k-j} -
\sum_{0\le j\le \frac{k}{2}} b_ja_{k-j}\Big|^p
\\ & \le 
\sum_{n=0}^\infty 2^n \sup_{2^n\le k<2^{n+1}} 
\Big(\big|(a*b)_k\big|+\Big|\sum_{0\le j\le \frac{k}{2}} b_ja_{k-j}\Big|\Big)^p
\\ & \le \nonumber
\sum_{n=0}^\infty 2^n \sup_{2^n\le k<2^{n+1}} 
\Big(\big|(a*b)_k\big|+
\Big(\sum_{0\le j\le \frac{k}{2}} |b_j|\Big)\sup_{0\le j\le \frac{k}{2}} |a_{k-j}|
\Big)^p
\\ & \le 
\sum_{n=0}^\infty 2^n \sup_{2^n\le k<2^{n+1}} 
\Big(\big|(a*b)_k\big|+
\|b\|_1\sup_{\frac{k}{2}\le j\le k} |a_{j}|\Big)^p
\\ & \le 
\sum_{n=0}^\infty 2^n  
\Big(\sup_{2^n\le k<2^{n+1}}\big|(a*b)_k\big|+
\|b\|_1\sup_{2^n\le k<2^{n+1}}\sup_{\frac{k}{2}\le j\le k} |a_{j}|\Big)^p
\\ & = 
\sum_{n=0}^\infty 2^n  
\Big(\sup_{2^n\le k<2^{n+1}}\big|(a*b)_k\big|+
\|b\|_1\sup_{2^{n-1}\le k<2^{n+1}} |a_{j}|\Big)^p.
\end{align*}

Minkowski's inequality and \eqref{eq-4-6} yield
\begin{align*} 
\Big(\sum_{n=0}^\infty 2^n \sup_{2^n\le k<2^{n+1}} 
\Big|\sum_{\frac{k}{2}<j\le k}b_ja_{k-j}\Big|^p\Big)^{1/p}
& \le 
\Big(\sum_{n=0}^\infty 2^n  
\sup_{2^n\le k<2^{n+1}}\big|(a*b)_k\big|^p\Big)^{1/p}
\\ & 
\quad +\|b\|_1\Big(\sum_{n=0}^\infty 2^n
\sup_{2^{n-1}\le k<2^{n+1}}|a_{j}|^p\Big)^{1/p}
\\ & \le
\|a*b\|_{\dpe}+ 3\|b\|_1\|a\|_{\dpe}.
\end{align*}
So,  \eqref{eq-4-1} holds.
\end{proof}


The equivalent norms for $\dpe$ given in \eqref{eq-2-3} 
and \eqref{eq-2-4} suggest, for each $1\le p<\infty$, to introduce the sequence space
\begin{equation}\label{eq-4-8}
d_{pp}:=\Big\{a=(a_n)_{n=0}^\infty\in\cno:
\sum_{n=0}^\infty 2^{np}  \sup_{2^{n}\le k< 2^{n+1}}|a_k|^p<\infty\Big\},
\end{equation}
equipped with the norm
\begin{align}\label{eq-4-9}
\|a\|_{\dpp} :=
\bigg(\sup_{k\ge0} |a_k|^p+\sum_{n=0}^\infty 2^{np} 
\sup_{2^n\le k<2^{n+1}} |a_k|^p\bigg)^{1/p},
\quad a\in\dpp.
\end{align}
The canonical vectors  $\{e_n:n\in\no\}$ form an unconditional basis in 
$\dpp$. To see this fix 
$a=(a_n)\subi\in \dpp$. For each $N\in\no$ let $n_0\in\no$ 
satisfy $2^{n_0}\le N<2^{n_0+1}$. Then,
for $N\to\infty$, we have
\begin{align*}
\Big\|a-\sum_{n=0}^N a_ne_n\Big\|^p_{\dpp} \le
\sup_{k>N} |a_k|^p 
+\sum_{n>n_0}^\infty 2^{np} \sup_{2^n\le k<2^{n+1}} |a_k|^p
\to 0.
\end{align*}
The bounded multiplier test ensures the unconditionality of the basis.


\begin{theorem}\label{t-4-2}
Let $p\in[1,\infty)$.  Then $\dpp\cap\ell^1\subsetneq\mdp$
with a continuous inclusion.
\end{theorem}

\begin{proof}
Since  $\mdone=d_1=d_{11}$, we only need to consider the case when $p\in(1,\infty)$. 
Fix $b\in\dpp\cap\ell^1$.
We apply Theorem \ref{t-4-1} by verifying that \eqref{eq-4-1} holds. 
Given $a\in\dpe$ we have 
\begin{align*} 
\Big|\sum_{\frac{k}{2}<j\le k}b_ja_{k-j}\Big|
& \le
\Big(\sum_{j=0}^{k/2}|a_{j}|\Big)\sup_{\frac{k}{2}<j\le k}|b_j|,
\quad k\ge1,
\end{align*}
and so  H\"older's inequality together with $\dpe\subseteq\ell^p$ yields
\begin{align*} 
\sup_{2^n\le k<2^{n+1}} 
\Big|\sum_{\frac{k}{2}<j\le k}b_ja_{k-j}\Big|^p
&\le 
\sup_{2^n\le k<2^{n+1}} 
\Big(\sum_{j=0}^{k/2}|a_{j}|\Big)^p\sup_{\frac{k}{2}<j\le k}|b_j|^p
\\ &\le 
\Big(\sum_{j=0}^{2^{n}-1}|a_{j}|\Big)^p\sup_{2^{n-1}<j< 2^{n+1}}|b_j|^p
\\ &\le 
2^{n(p/q)}\|a\|_{\dpe}^p\sup_{2^{n-1}\le j< 2^{n+1}}|b_j|^p.
\end{align*}
Hence, arguing as in \eqref{eq-4-6}, it follows that 
\begin{align}\label{eq-4-10}
\sum_{n=0}^\infty 2^n \sup_{2^n\le k<2^{n+1}} 
\Big|\sum_{\frac{k}{2}<j\le k}b_ja_{k-j}\Big|^p
&\le 
\sum_{n=0}^\infty 2^n 2^{n(p/q)}\|a\|_{\dpe}^p\sup_{2^{n-1}\le j< 2^{n+1}}|b_j|^p
\\ & =\nonumber
\|a\|_{\dpe}^p \sum_{n=0}^\infty 2^{np}\sup_{2^{n-1}\le j< 2^{n+1}}|b_j|^p
\\ & \le
3\|a\|_{\dpe}^p \sum_{n=0}^\infty 2^{np}\sup_{2^{n}\le j< 2^{n+1}}|b_j|^p<\infty,\nonumber
\end{align}
which is finite since $b\in\dpp$. So, $\dpp\cap\ell^1\subseteq\mdp$.

In view of \eqref{eq-4-5} and \eqref{eq-4-9}, it follows from \eqref{eq-4-7}
and \eqref{eq-4-10} that there exists a constant  $K>0$ such that
$$
\|b*a\|_{\dpe}\le K\|a\|_{\dpe}\max\big\{\|b\|_1,\|b\|_{\dpp}\big\},
\quad a\in\dpe.
$$
Since the space $\dpp\cap\ell^1$ is normed by 
$\|b\|_{\dpp\cap\ell^1}:=\max\{\|b\|_1,\|b\|_{\dpp}\}$, it follows 
that the natural inclusion $\dpp\cap\ell^1\subseteq\mdp$ is continuous.

It remains to show that there exists $b\in\mdp\setminus\dpp$. Consider $b=(b_n)\subi$ defined by
$b_n=1/n$ for $n=2^k$ with $k\in\no$, and $b_n=0$ elsewhere. Then $b\not\in\dpp$ since
$$
\sum\subi 2^{np} \sup_{2^n\le k<2^{n+1}} |b_k|^p=\sum\subi 2^{np} \Big(\frac{1}{2^n}\Big)^p=\infty.
$$
However, $b\in\mdp$. Indeed, via Theorem \ref{t-4-1}
and the fact that $b\in\ell^1$ we have
$$
\sum_{n=0}^\infty 2^n \sup_{2^n\le k<2^{n+1}} 
\Big|\sum_{\frac{k}{2}<j\le k}b_ja_{k-j}\Big|^p
= \sum_{n=0}^\infty 2^n \Big|\frac{a_0}{2^n}\Big|^p
<\infty,\quad a\in\dpe.
$$
\end{proof}


The containment  $d_1\subseteq \dpp$ follows directly from \eqref{eq-4-8}
because of \eqref{eq-2-3}, \eqref{eq-4-9} and
$$
\sum_{n=0}^\infty 2^{np}  \sup_{2^{n}\le j< 2^{n+1}}|a_j|^p
\le \Big(\sum_{n=0}^\infty 2^{n}  \sup_{2^{n}\le j< 2^{n+1}}|a_j|\Big)^p.
$$
Thus, Theorem \ref{t-4-1} and the fact that $d_1=\mdone$ imply the following result
(a strengthening of part of Theorem \ref{t-3-8}).

\begin{corollary}\label{c-4-3}
Let $p\in[1,\infty)$.   The following continuous inclusion holds:
$$
d_1\subseteq\mdp.
$$
\end{corollary}


Let $\hdd$ denote the algebra, under pointwise multiplication, of all $\C$-valued 
functions which are holomorphic in some open set containing $\DD$.

\begin{corollary}\label{c-4-4}
Let $p\in[1,\infty)$. The following  inclusions hold:
$$
\Big\{b=(b_n)\subi: f_b\in \hdd\Big\}\subseteq d_1\subseteq\mdp.
$$
\end{corollary}

\begin{proof}
Given $f_b\in\hdd$, the power series  of $f_b$ 
has radius of convergence $r>1$ and so its Taylor coefficients 
satisfy $|b_n|\le C/r^n$, for some $C>0$ and all $n\in\no$. Hence, 
$b\in d_1\subseteq\mdp$ for all $p\in[1,\infty)$.
\end{proof}


\begin{corollary}\label{c-4-5} 
Let $p\in[1,\infty)$. For $b=(b_n)\subi$ belonging to any one of the
spaces $\lwp$ or $\dpp\cap\ell^1$ or $d_1$, it is the case, for $N\to\infty$, that
$$
\Big\|b-\sum_{n=0}^N b_ne_n\Big\|_{\mdp}\to0.
$$
Equivalently, 
$$
\Big\|T_b-\sum_{n=0}^N b_nS^n\Big\|_{\modp}\to0.
$$
\end{corollary}

\begin{proof}
The sequence $\{e_n:n\in\no\}$ is a basis for each of these spaces. This,
together with Proposition \ref{p-3-7}, Theorem \ref{t-4-2} and Corollary \ref{c-4-3},
proves the result.
\end{proof}


\begin{remark}
We compare the various subspaces of $\mdp$   which have already
appeared.

(i) For every $p\in[1,\infty)$ the  spaces
$d_1$ and $\lwp$ are different. Indeed,  $b=(b_n)\subi$ given by 
$b_n:=1/(n+1)^{1+\frac1p}$, for $n\in\no$, satisfies
$b\in d_1$ but $b\not\in\lwp$. So, $b\in\mdp\setminus\lwp$. On the other hand, 
the example $b$ in the proof of Theorem \ref{t-4-2} satisfies $b\in\lwp$ but $b\not\in d_1$ as $b\not\in\dpp$. 
So, $b\in\mdp\setminus d_1$.

(ii) For every $p\in(1,\infty)$ we have $\dpp\subsetneq \dpe$. The containment is direct from 
\eqref{eq-2-3} and \eqref{eq-4-8}. To see that
it is strict, consider again the example $b$ in the proof of Theorem \ref{t-4-2}.
Then $b\in\dpe$ but $b\not\in \dpp$.

(iii)  For every $p\in(1,\infty)$ we have $\ell^1\not\subseteq \dpp$. The 
proof of Corollary \ref{c-3-6}(ii) yields $\ell^1\not\subseteq \dpe$.
To see that $\dpp\not\subseteq\ell^1$,
consider $b=(b_n)\subi$ with $b_0=0$ and $b_n=1/(k2^{k})$ when $2^k\le n<2^{k+1}$ and $k\in\no$.
Then $b\in\dpp$ but $b\not\in \ell^1$. This sequence $b$ shows that 
$d_1\subsetneq \dpp\cap\ell^1$, since it satisfies
$b\in\dpp\cap\ell^1$ and $b\not\in d_1$.

(iv) For every $p\in[1,\infty)$ the  spaces
$\dpp$ and $\lwp$ are different. Indeed,  $b=(b_n)\subi$ given by 
$b_n:=1/(n+1)^{1+\frac1p}$, for $n\in\no$, satisfies
$b\in \dpp$ but $b\not\in\lwp$. 
On the other hand, 
the example $b$ in the proof of Theorem \ref{t-4-2}  satisfies $b\in\lwp$ and  $b\not\in \dpp$. 
\end{remark}


\section{Spectral  properties of $\mdp$}
\label{S6}

It was noted in Section \ref{S1} that the multiplier algebra 
$\mathscr{M}(ces_p)=\ell^1$ for every $1<p<\infty$.
For elements $b\in\ell^1$, the spectrum of the corresponding operator
$T_b\in\mathscr{L}(ces_p)$ is precisely known,
\cite[Theorem 2]{ricker}. The proof requires a 
knowledge of the spectrum of the right-shift $S\in\mathscr{L}(ces_p)$,
which is identified in \cite[Proposition 6]{ricker}.
The aim of this section is to investigate the spectrum of
multiplier operators $T_b\in\mdp$ for $1\le p<\infty$.
Due to the more involved nature of the Banach algebras $\mdp$ 
this is significantly more complicated than the situation 
for $ces_p$. We begin with the right-shift $S\in\ldp$.
The spectrum of $S\in{\mathscr{L}(d_p)}$ is well known,
\cite[VII Proposition 6.5]{conway}.


\begin{proposition}\label{p-5-1}
Let $p\in[1,\infty)$. The right-shift operator $S\colon \dpe\to\dpe$ satisfies
\begin{equation}\label{eq-5-1}
\sigma(S;\ldp)=\overline{\D}.
\end{equation}
Moreover, the point spectrum
\begin{equation*}
\sigma_{\text{pt}}(S;\ldp)=\emptyset
\end{equation*}
and  the residual spectrum satisfies
\begin{equation*}
\D\subseteq \sigma_{\text{r}}(S;\ldp).
\end{equation*}
Whenever $p\in(1,\infty)$,   the continuous spectrum satisfies
\begin{equation}\label{eq-5-2}
\sigma_{\text{c}}(S;\ldp)=\overline{\D}\setminus\D.
\end{equation}
\end{proposition}

\begin{proof}
The proof proceeds via a series of steps.  
All steps, but for for the last one, concern $p\in[1,\infty)$.

\textit{Step 1.}   
We have that
$$
\sigma_{\text{pt}}(S;\ldp)=\emptyset.
$$
To prove this, suppose that $\lambda\in\sigma_{\text{pt}}(S;\ldp)$. Then there exist
$0\not=a\in\dpe$ such that $Sa=\lambda a$. Since $a\in\ell^p$ this implies that
$a$ is an eigenvalue of  $S\colon\ell^p\to \ell^p$. This cannot be
since $\sigma_{\text{pt}}(S;\mathscr{L}(\ell^p))=\emptyset$; see
\cite[Proposition VII.6.5]{conway}.

\textit{Step 2.}   
For the range $R(S-\lambda I)$ of $S-\lambda I$ it is the case that 
$$
e_0\not\in R(S-\lambda I)\subseteq\dpe,\quad \lambda\in\overline{\D}.
$$
To prove this, fix $\lambda\in\overline{\D}$. Suppose there exists $a\in\dpe$ such that 
$(S-\lambda I)a=e_0$. Necessarily $a\not=0$.  If $\lambda=0$, then $Sa=e_0$, which is
impossible. 
For $0<|\lambda|\le1$ we have 
$$
-\lambda a_0=1,\quad -\lambda a_{n+1}=a_n,\quad n\in\no.
$$
Proceeding recursively yields $a_n=1/\lambda^{n+1}$ for $n\in\no$. But, then $a\not\in\dpe$ as
$1/|\lambda|\ge1$.

\textit{Step 3.}  The same calculations as in Step 2, for $\ell^p$ in place of $\dpe$ 
and the right-shift operator $S\in\mathscr{L}(\ell^p)$ show that
$$
e_0\not\in R(S-\lambda I)\subseteq\ell^p,\;\;\lambda\in\DD.
$$

\textit{Step 4.}   
For each $\lambda\in\D$, it is the case that
$$
e_0\not\in \overline{R(S-\lambda I)}\subseteq\dpe,
$$
where the bar denotes closure.
To prove this, fix $\lambda\in\D$. Suppose, on the contrary, that there exists a
sequence $\{a^m\}_{m=0}^\infty\subseteq \dpe$ such that  $(S-\lambda I)a^m\to e_0$ in $\dpe$.
Then also $e_0\in\ell^p$ and the sequence $\{a^m\}_{m=0}^\infty\subseteq \ell^p$ satisfies  
$(S-\lambda I)a^m\to e_0$ in $\ell^p$. But, the range
$R(S-\lambda I)$ is \textit{closed} in $\ell^p$; see Proposition VII.6.5 in \cite{conway}.
Hence, $e_0\in R(S-\lambda I)\subseteq\ell^p$ which contradicts Step 3.

\textit{Step 5.}  
For the residual spectrum we have the inclusion
$$
\D\subseteq \sigma_{\text{r}}(S;\ldp).
$$
To prove this note,  by Step 1,  that $S-\lambda I$ is injective for every $\lambda\in\D$.
Accordingly, for each $\lambda\in\D$, Step 4 shows that $\overline{R(S-\lambda I)}\not=\dpe$
and hence, that $\lambda\in \sigma_{\text{r}}(S;\ldp)$.

\textit{Step 6.}  
The claim is   that
$$
\sigma(S;\ldp)\subseteq \overline{\D}.
$$
To prove this, recall that $\|S^n\|_{\ldp}=(n+1)^{1/p}$ for $n\in\no$.
Accordingly, the spectral radius 
$r(S)=\lim_{n}\|S^n\|_{\ldp}^{1/n}=1$
from which the result follows, \cite[I Theorem 5.8]{bonsall-duncan}.

\textit{Step 7.} 
The identity \eqref{eq-5-1} is valid, that is,
$$
\sigma(S;\ldp)= \overline{\D}.
$$
To prove this, note that Steps 5 and 6 yield
$$
\D\subseteq \sigma_{\text{r}}(S;\ldp) \subseteq \sigma(S;\ldp)\subseteq\overline{\D}.
$$
Since the spectrum of $S$ is a closed set in $\C$ the desired conclusion follows.

\textit{Step 8.}  For  
every $\lambda\in\C\setminus\{0\}$ it is the case that
\begin{equation*}\label{eq-5-3}
\big\{-\lambda e_0+\frac{1}{\lambda^n} e_{n+1}:n\in\no\big\}\subseteq R(S-\lambda I)
\subseteq\dpe.
\end{equation*}
To verify this define, for each $n\in\no$, the element
$$
a^{[n]}:=\sum_{j=0}^n \frac{1}{\lambda^j} e_{j}
=\Big(1,\frac1\lambda,\dots,\overbrace{\frac{1}{\lambda^n}}^{\text{position } n+1},
0,\dots\Big)\in\dpe.
$$
Direct calculation yields
$$
(S-\lambda I)a^{[n]}=
\Big(-\lambda,0,\dots,0,\overbrace{\frac{1}{\lambda^{n}}}^{\text{position } n+2},0,\dots\Big)
=-\lambda e_0+\frac{1}{\lambda^n} e_{n+1}.
$$

\textit{Step 9.}   Consider now $p\in(1,\infty)$. Then 
$$
\sigma_{\text{c}}(S;\ldp) = \overline{\D}\setminus \D.
$$
To prove this, recall that $\dpe^*=ces_q$, with $\frac1p+\frac1q=1$. 
Fix $\lambda\in \overline{\D}\setminus \D$.
Let $y^*=(y_n)\subi\in\dpe^*$ satisfy
\begin{equation}\label{eq-5-3}
\big\langle-\lambda e_0+\frac{1}{\lambda^n} e_{n+1},y^*\big\rangle=0,\quad n\in\no.
\end{equation}
Substituting $n=0,1,\dots$ successively into \eqref{eq-5-3} 
yields $y_n=\lambda^n y_0$, for all $n\in\no$, and so
$y^*=(y_0\lambda^n)\subi$. Then $|y^*|=(|y_0|)\subi\in\dpe^*=ces_q$. 
The definition of $ces_q$
in \eqref{eq-2-5} implies that $|y^*|=\ce|y^*|\in \ell^q$ which implies 
that $y_0=0$, that is, $y^*=0$.

Now let $y^*\in\dpe^*$ satisfy $\langle a,y^*\rangle=0$ for all $a\in R(S-\lambda I)$. According
to Step 8, $y^*$  also satisfies \eqref{eq-5-3} and hence, $y^*=0$. It follows that
$\overline{R(S-\lambda I)}=\dpe$. Since $\lambda\in \sigma(S; \ldp)$, due to Step 7, and
$S-\lambda I$ is injective (see Step 1), we can conclude that 
$\lambda\in \sigma_{\text{c}}(S;\ldp)$. That is, 
$\overline{\D}\setminus \D\subseteq \sigma_{\text{c}}(S;\ldp)$. Now Steps 5 and 7 yield
$\sigma_{\text{c}}(S;\ldp) = \overline{\D}\setminus \D$.

The proof is thereby complete.
\end{proof}


The omission of $p=1$ in \eqref{eq-5-2} is necessary, as seen by the following result.

\begin{proposition}
For $p=1$ we have that
\begin{equation*}
\sigma(S;\mathscr{L}(d_1))=\sigma_{\text{r}}(S;\mathscr{L}(d_1))=\overline{\D}.
\end{equation*}
In particular,
\begin{equation*}
\sigma_{\text{pt}}(S;\mathscr{L}(d_1))=\sigma_{\text{c}}(S;\mathscr{L}(d_1))=\emptyset.
\end{equation*}
\end{proposition}

\begin{proof}
According to Proposition \ref{p-5-1} we only need to show that if $|\lambda|=1$,
then $\lambda\in\sigma_{\text{r}}(S;\mathscr{L}(d_1))$.  Recall that $d_1^*=(ces_0)^{**}=ces_\infty$,
\cite[Remark 6.3]{curbera-ricker-ieot-2014}. Set $y^*:=(\lambda^n)\subi$. Observe that $|y^*|=(1)\subi$ and,
for $\ce$ the \ces averaging operator, that $\ce|y^*|=(1)\subi\in\ell^\infty$. Hence, by definition $y^*\in ces_\infty=d_1^*$.

Let $a\in d_1$ be arbitrary. Then
\begin{align*}
\langle(S-\lambda I)a,y^*\rangle &=
\big\langle(-\lambda a_0,a_0-\lambda a_1, a_1-\lambda a_2,\dots),(1,\lambda,\lambda^2,\dots)\big\rangle
\\ & =
-\lambda a_0+\lambda(a_0-\lambda a_1)+\lambda^2(a_1-\lambda a_2)+\cdots 
\\ &=0.
\end{align*}
That is, $y^*\not=0$ in $d_1^*$ satisfies $\langle u,y^*\rangle=0$ 
for all $u\in R(S-\lambda I)\subseteq d_1$.
Accordingly, $\overline{R(S-\lambda I)}\not=d_1$. Since $S-\lambda I$ is  injective,
we can conclude that $\lambda\in\sigma_{\text{r}}(S;\mathscr{L}(d_1))$.
\end{proof}


The above knowledge of the spectrum for the right-shift operator has implications 
for other multipliers. Given $f\in\hdd$, let $b_f=(b_n)\subi$ denote the sequence of its
Taylor coefficients.

\begin{proposition}\label{p-5-3}
Let $p\in[1,\infty)$. For every $f\in\hdd$ we have that $b_f\in\mdp$ and
\begin{equation*}
\sigma(T_{b_f};\modp)=\sigma(T_{b_f};\ldp)=f(\DD).
\end{equation*}
\end{proposition}

\begin{proof} 
Fix $f\in\hdd$. We know   (cf.\ Corollary \ref{c-4-4}) that 
$b_f\in\mdp$ and so $T_{b_f}\in\modp$.
Via the functional calculus for unital Banach algebras, 
\cite[Ch.I, \S7]{bonsall-duncan}, \cite[Ch.10 \& 11]{rudin}, the operator 
$f(S)\in\modp$ is defined by the Cauchy integral formula
$$
f(S):=\frac{1}{2\pi i}\int_\gamma f(z)(zI-S)^{-1}dz
$$ 
for a suitable contour $\gamma$ surrounding $\DD=\sigma(S; \modp)$,
where we use Remark \ref{r-3-2}(ii) and \eqref{eq-5-1}.

Fix $n\in\no$. Given $z\in\gamma$ a direct calculation yields (as $|z|>1$) that 
$$
(zI-S)^{-1}e_n=\Big(0,\dots,0,\overbrace{\frac{1}{z}}^{\text{position } n}, 
\frac{1}{z^2},\frac{1}{z^3},0,\dots\Big)
\in d_1\subseteq \dpe.
$$
Accordingly,
$$
f(S)e_n=\frac{1}{2\pi i}\int_\gamma f(z)(zI-S)^{-1}e_n\,dz
=\sum_{k=0}^\infty \frac{1}{2\pi i}\int_\gamma \frac{f(z)}{z^{k+1}}\,dz\cdot e_{k+n}.
$$
Since $b_f=(\frac{1}{2\pi i}\int_\gamma \frac{f(z)}{z^{k+1}}\,dz)_{k=0}^\infty$, it follows that
$$
f(S)e_n=\Big(0,\dots,0,\overbrace{b_0}^{\text{position } n}, b_1,b_2,\dots\Big)=b_f*e_n.
$$
But, $b_f\in\mdp$, that is, $T_{b_f}\in\modp$ and so $f(S)e_n=T_{b_f}e_n$ for all $n\in\no$.
Since $\{e_n:n\in\no\}$ is basis for $\dpe$, we can 
conclude that $f(S)=T_{b_f}$. By the spectral mapping theorem
for $f(S)$ and \eqref{eq-5-1} we have
$$
\sigma(f(S);\ldp)=f(\sigma(S;\ldp))=f(\DD).  
$$ 
Since $\sigma(f(S);\ldp)=\sigma(f(S);\modp)=\sigma(T_{b_f};\modp)$,
the proof is complete.
\end{proof}


\begin{proposition}\label{p-5-4}
The maximal ideal space of $\mathscr{M}(d_1)$ is homeomorphic to $\DD$.
Moreover, for each $b\in \mathscr{M}(d_1)=d_1$, its spectrum is given by
\begin{equation*}
\sigma(b;\mathscr{M}(d_1))=\sigma(T_b;\mathscr{M}_{\text{op}}(d_1))=f_b(\DD). 
\end{equation*}
\end{proposition}

\begin{proof}
Recall that $d_1$ is an algebra, that is, $\mathscr{M}(d_1)=d_1$ with equivalence of norms. 
Moreover, the unital Banach algebra $\mdone$ is generated by $e_1$. To see this,
let $b=(b_n)\subi\in\mdone=d_1$. 
Recall that $e_m=e^m_1$ for all $m\ge1$ and so each element 
 $b^n:=b_0e_0+\sum_{j=1}^nb_je_j$, for $n\in\no$, belongs
to the algebra $\langle e_0,e_1\rangle$ generated by $e_0$ and $e_1$. Since $\{e_n:n\in\no\}$ is a basis 
for $d_1$ and $\mdone=d_1$, it follows that $b^n\to b$ in the norm of $d_1$ and hence, in the
norm of $\mdone$. So, the closure of $\langle e_0,e_1\rangle$ in $\mdone$ is $\mdone$.

Theorem 2 on p. 98 of \cite{bonsall-duncan} implies that the maximal ideal space $\Phi$ of
$\mdone$ is homeomorphic with the spectrum $\sigma(e_1;\mdone)$ of the generator $e_1$.
Since $\mdone$ is isometric  to $\mathscr{M}_{\text{op}}(d_1)$ we know from 
Proposition \ref{p-5-1} that
\begin{equation*}
\sigma(e_1;\mathscr{M}(d_1))=\sigma(T_{e_1};\mathscr{M}_{\text{op}}(d_1))=
\sigma(S;\mathscr{M}_{\text{op}}(d_1))=\sigma(S;\mathscr{L}(d_1))=\DD.
\end{equation*}
More explicitly, each $z\in\DD\simeq \Phi$ defines the multiplicative, linear functional
on $\mdone$ via point evaluation, namely
$$
b\mapsto f_b(z),\quad b\in\mdone=d_1.
$$
Since $b\in d_1\subseteq\ell^1$, the continuity is immediate from
$|f_b(z)|=|\sum_{n=0}^\infty b_n z^n|\le \sum_{n=0}^\infty |b_n|\le\|b\|_{d_1}$, for
$b\in \mdone$. The Gelfand transform $\hat b\colon\Phi\to\C$, of each $b\in\mdone$ is given
by $\hat b(z)=f_b(z)$, for $z\in\DD$. It follows from Theorem 11.9.(c) in \cite{rudin}
that $\sigma(b;\mathscr{M}(d_1))=\hat b(\Phi)=f_b(\DD)$ for each $b\in\mdone$.
\end{proof}


Fix $p\in[1,\infty)$ and let $\bodp$ denote the closure in $\modp$ of the
algebra $\langle I,S\rangle$ consisting of all operators which are 
polynomials in $S$.


\begin{proposition}\label{p-5-5} 
Let $p\in[1,\infty)$. The maximal ideal space of $\bodp$ is homeomorphic to $\DD$.
Moreover, for each $T_b\in \bodp$, that is, for each $b\in\mdp$  such that $T_b\in\bodp$,
its spectrum is given by
\begin{equation*}
\sigma(T_b;\bodp)=f_b(\DD). 
\end{equation*}
\end{proposition}

\begin{proof}
The discussion at the beginning of the proof of Proposition \ref{p-5-4} shows that 
$\mathscr{A}(S,d_1)
=\mathscr{M}_{\text{op}}(d_1)=d_1$ 
and so Proposition \ref{p-5-4} establishes the desired identity.

Next  consider  $p\in(1,\infty)$. Since the multiplication in any Banach algebra 
is jointly continuous, it follows that $\bodp$ is a \textit{closed} subalgebra of
$\modp$. Moreover, $\sigma(S;\modp)=\DD$; see Remark \ref{r-3-2}(ii) and 
Proposition \ref{p-5-1}. Since $\C\setminus\DD$ is a connected set, it follows 
from \cite[I Proposition 5.14]{bonsall-duncan} that also $\sigma(S;\bodp)=\DD$.
In particular, the maximal ideal space of $\bodp$ is homeomorphic to $\DD$
(cf.\   \cite[II Theorem 19.2]{bonsall-duncan}) and so, for any polynomial $f$, we have that
$$
\sigma(f(S);\bodp)=\sigma(f(S);\modp)=f(\DD).
$$

Every $T\in\bodp\subseteq\modp$ is of the form $T=T_b$ for some unique element $b\in\mdp$. 
Each $z\in\DD$  defines the linear, multiplicative functional on $\bodp$ via
$$
T_b\mapsto f_b(z),\quad T_b\in\bodp,
$$
which is automatically continuous, \cite[II Proposition 16.3]{bonsall-duncan}.
The Gelfand transform $\widehat{T_b}\colon\DD\to\C$, of each $T_b\in\bodp$,
is given by $\widehat{T_b}(z)=f_b(z)$, for $z\in\DD$. Again by Theorem 11.9(c)
in \cite{rudin} we can conclude that $\sigma(T_b;\bodp)=\widehat{T_b}(\DD)$.
\end{proof}


\begin{remark}\label{r-6-6}
(i) Let $b\in\mdone$ belong to the radical. Proposition \ref{p-5-4} together with 
Theorem 11.9.(c) in \cite{rudin} imply, for the Gelfand transform $\hat b$, 
that $\|\hat b\|_\infty=0$, that is, $f_b(\DD)=0$ and so $b=0$. Hence,
$\text{rad}(\mdone)=\{0\}$, that is,
$\mdone$ is \textit{semisimple}. An analogous argument (now using 
Proposition \ref{p-5-5}) shows   that also $\bodp$ is a 
semisimple algebra for all $p\in[1,\infty)$.

(ii) Given $p\in[1,\infty)$, which elements $b\in\mdp$ satisfy $T_b\in\bodp$?
According to Corollary \ref{c-4-5}, this includes the space $d_1$ (hence,
also the Taylor coefficients $b_f$ of any function $f\in\hdd$ via Corollary \ref{c-4-4}),
the weighted space $\lwp$ and $\dpp\cap\ell^1$. Actually, for every
$b=(b_n)\subi$ belonging to any one of these spaces, the approximation of $T_b$ can be achieved by
using the \textit{Taylor polynomials} of $b$. That is, for $n\to\infty$, we have
$$
\Big\|T_b-\sum_{j=0}^n b_jS^j\Big\|_{\bodp}
=\Big\|T_b-\sum_{j=0}^n b_jS^j\Big\|_{\modp}\to0.
$$
\end{remark}


The following identities occur in Proposition \ref{p-5-3}, namely
$$
\sigma(T_{b_f};\bodp)=\sigma(T_{b_f};\modp)=f(\DD),\quad f\in\hdd.
$$
For certain other multipliers   an inclusion is possible.

\begin{proposition}\label{p-5-7}
Let $p\in[1,\infty)$ and $b\in \mdp$ satisfy
\begin{equation}\label{eq-5-4}
\Big\|T_b-\sum_{j=0}^n b_jS^j\Big\|_{\modp}\to0\quad \text{for $n\to\infty$}.
\end{equation}
Then
$$
\sigma(T_b;\bodp)=\Big\{\sum\subi b_n\lambda^n:\lambda\in\DD\Big\}
\subseteq \sigma(T_b;\modp)=\sigma(T_b;\ldp).
$$
\end{proposition}

\begin{proof}
Fix $\lambda\in\DD$. Since $b\in\ell^1$ (cf.\ Theorem \ref{t-3-5}) the series
$\sum_{j=0}^\infty b_j\lambda^j$ converges absolutely in $\C$. Define
$\alpha_n:=\sum_{j=0}^nb_j\lambda^j$, for $n\in\no$, in which case 
$\alpha_n\to\alpha:=\sum_{j=0}^\infty b_j\lambda^j$ for $n\to\infty$.
Moreover, setting $R_n:=\sum_{j=0}^nb_jS^j$ we have that $R_n\in\modp$
and so
$$
\sigma(R_n;\modp)=\Big\{\sum_{j=0}^n b_j z^j:z\in\DD=\sigma(S;\modp)\Big\},
\quad n\in\no.
$$
That is, $\alpha_n\in \sigma(R_n;\modp)$ for $n\in\no$. For $\mathscr{A}:=\modp$
it follows from \eqref{eq-5-4} that $R_n\to T_b$ in $\mathscr{A}$ and 
so  \cite[Ex. 5, p.199]{conway} implies that
$\sum_{j=0}^\infty b_j\lambda^j\in \sigma(T_b;\modp)$.
\end{proof}




\begin{thebibliography}{99}


\bibitem{belinskii}
Belinskii, E. S., Liflyand, E. R., Trigub, R. M.,
\emph{The Banach algebra $A^*$ and its properties},
J. Fourier Anal. Appl.,   \textbf{3} (1997), 103--129. 

\bibitem{bennett}
Bennett, G.,
\emph{Factorizing the classical inequalities},
Mem. Amer. Math. Soc.,  \textbf{120} (576), (1996), 1--130. 

\bibitem{beurling}
Beurling, A., \textit{On the spectral synthesis of bounded functions}, 
Acta  Math.,  \textbf{81} (1949), 225--238.


\bibitem{bottcher-silbermann}
B\"ottcher, A., Silbermann, B.,
\emph{Analysis of Toeplitz Operators},
Springer,  Berlin Heidelberg, 2006. 

\bibitem{bonet-ricker}
Bonet, J., Ricker, W. J.,
\emph{Operators acting in the dual spaces of discrete
Ces\`{a}ro  spaces},
Monatsh. Math.,  \textbf{191} (2020), 487--512. 

\bibitem{bonsall-duncan}
Bonsall, F. F., Duncan, J.,
\emph{Complete Normed Algebras},
Springer,  Berlin Heidelberg, 1973. 

\bibitem{calderon}
Calder\'on, A. P., 
\emph{Intermediate spaces and interpolation, the complex method}, 
Studia Math. \textbf{24} (1964), 113--190.


\bibitem{cheng-etal}
Cheng, R., Mashreghi, J., Ross, W.T., 
\emph{Function Theory and $\ell^p$ Spaces},
Univ. Lecture Series \textbf{75}, 
Amer. Math. Soc., Providence, R.I., 2020.

\bibitem{conway}
Conway, J. B., 
\emph{A Course in Functional Analysis}, 2nd Ed.,
Springer, New York, 1990. 

\bibitem{curbera-ricker-ieot-2014}
Curbera, G. P., Ricker, W. J., 
\emph{Solid extensions of the Ces\`{a}ro  operator on $\ell^p$ and $c_0$},
Integr. Equ. Operator Theory, \textbf{80} (2014), 61--77.

\bibitem{curbera-ricker-jmaa-2022}
Curbera, G. P., Ricker, W. J.,
\emph{Fine spectra and compactness 
of generalized Ces\`{a}ro  operators in Banach lattices in $\cno$},
J. Math. Anal. Appl., \textbf{507} (2022), 125854, 31 pp. 

\bibitem{grosse-erdmann}
Grosse-Erdmann, K.-G., 
\emph{The Blocking Technique, Weighted Mean Operators and Hardy's 
Inequality}, Lecture Notes in Mathematics \textbf{1679}, Springer,
Berlin-Heidelberg, 1998. 

\bibitem{lesnik-maligranda}
Lesnik, K., Maligranda, L., \emph{Interpolation of Ces\`{a}ro, Copson and Tandori spaces},
 Indag. Math., \textbf{27} (2016), 764--785.


\bibitem{maligranda}
Maligranda, L., \emph{Orlicz Spaces and Interpolation},
Seminars in Mathematics \textbf{5}, Univ. Estadual de Campinas, 
Campinas SP, Brazil 1989.

\bibitem{naimark}
Naimark, M. A., 
\emph{Normed Algebras}, 
Wolters-Noordhoff Publishing, Groningen, 1972. 

\bibitem{nikolskii}
Nikolskii, N. K.,
\emph{Spaces and algebras of Toeplitz matrices operating in $\ell^p$},  
Siber. Math. J.,  \textbf{7} (1966), 118--126. 

\bibitem{nursultanov}
Nursultanov, E., Tikhonov, S.,
\emph{Wiener-Beurling spaces and their properties},  
Bull. Sci. Math.,  \textbf{159} (2020), 102825, 20 pp.

\bibitem{ricker}
Ricker, W. J.,
\emph{Convolution operators acting in discrete
Ces\`{a}ro  spaces},
Arch. Math.,  \textbf{112} (2019), 71--82. 

\bibitem{rudin}
Rudin, W., 
\emph{Functional Analysis}, 2nd Ed.,
McGraw-Hill, Singapore, 1991. 

\bibitem{shestakov}
Shestakov,  V. A., \emph{Interpolation of linear operators in 
spaces of measurable functions}, 
Funct. Anal.  Appl., \textbf{8} (1974), 274--275. 


\end{thebibliography}
\end{document}